\input amstex 
\documentstyle{amsppt}
\input bull-ppt
\keyedby{bull479/kta}

\topmatter
\cvol{31}
\cvolyear{1994}
\cmonth{July}
\cyear{1994}
\cvolno{1}
\cpgs{1-14}
\title Selberg's Conjectures and Artin 
$L$-functions\endtitle
\author M. Ram Murty\endauthor
\shortauthor{M. Ram Murty}
\shorttitle{Selberg's conjectures and Artin $L$-functions}
\address Department of Mathematics, McGill University, 
Montreal,
H3A 2K6 Canada\endaddress
\ml murty\@math.mcgill.ca\endml
\date September 4, 1992, and, in revised form, September 
28, 1993\enddate
\tk Research partially supported by an NSERC grant\endtk
\endtopmatter

\document
\ah 1. Introduction\endah In its comprehensive form, an 
identity between an
automorphic $L$-function and a ``motivic'' $L$-function is 
called a
reciprocity law. The celebrated Artin reciprocity law is 
perhaps the
fundamental example. The conjecture of Shimura-Taniyama 
that every elliptic curve
over $Q$ is ``modular'' is certainly the most intriguing 
reciprocity
conjecture of our time. The ``Himalayan peaks'' that hold 
the secrets
of these nonabelian reciprocity laws challenge humanity, 
and, with the
visionary Langlands program, we have mapped out before us 
one means of
ascent to those lofty peaks. The recent work of Wiles 
suggests that an
important case (the semistable case) of the 
Shimura-Taniyama conjecture is on the
horizon and perhaps this is another means of ascent. In 
either case, a long
journey is predicted. To paraphrase the cartographer, it 
is not a journey
for the faint-hearted. Indeed, there is a forest to 
traverse, ``whose trees
will not fall with a few timid blows. We have to take up 
the double-bitted
axe and the cross-cut saw and hope that our muscles are 
equal to them.''

At the 1989 Amalfi meeting, Selberg [S] announced a series 
of conjectures
which looks like another approach to the summit. Alas, 
neither path seems
the easier climb. Selberg's conjectures concern Dirichlet 
series, which admit
analytic continuations, Euler products, and functional 
equations.

The Riemann zeta function is the simplest example of a 
function in the family
$\scr S$ of functions $F(s)$ of a complex variable $s$ 
satisfying the
following properties:

(i) (Dirichlet series) For $\roman{Re}(s)>1$, 
$F(s)=\sum^\infty
_{n=1}a_n/n^s$, where $a_1=1$, and we will write 
$a_n(F)=a_n$ for the
coefficients of the Dirichlet series.

(ii) (Analytic continuation) $F(s)$ extends to a 
meromorphic function so that
for some integer $m\geq 0$, $(s-1)^mF(s)$ is an entire 
function of finite
order.

(iii) (Functional equation) There are numbers $Q>0$, 
$\alpha_i>0$, $\roman
{Re}(r_i)\geq 0$, so that
$$\Phi(s)=Q^s\prod^d_{i=1}\Gamma(\alpha_is+r_i)F(s)$$
satisfies $\Phi(s)=w\overline{\Phi(1-\overline{s})}$ for 
some complex
number $w$ with $|w|=1$.

(iv) (Euler product) $F(s)=\prod_pF_p(s)$ where $F_p(s)=\exp
(\sum^\infty_{k=1}b_{p^k}/p^{ks})$ where 
$b_{p^k}=O(p^{k\theta})$ for some
$\theta<1/2$, where $p$ denotes a prime number (here and 
throughout this
paper).

(v) (Ramanujan hypothesis) $a_n=O(n^\epsilon)$ for any 
fixed $\varepsilon
>0$.

Note that the family $\scr S$ is multiplicatively closed 
and so is a
multiplicative monoid.

All known examples of elements in $\scr S$ are automorphic 
$L$-functions. In
all of these cases, $F_p(s)$ is an inverse of a polynomial 
in $p^{-s}$
of bounded degree.

Selberg [S] introduced this family to study the value 
distribution of
finite linear combinations of Dirichlet series with Euler 
products and
functional equations. For this purpose, he introduced the 
important concept
of a primitive function and made significant conjectures 
about them.

A function $F\in\scr S$ is called {\it primitive} if the 
equation
$F=F_1F_2$ with $F_1,F_2\in\scr S$ implies $F=F_1$ or 
$F=F_2$. As we shall
see below, one of the most serious consequences of the 
Selberg
conjectures is that $\scr S$ has unique factorization into 
primitive
elements. It is not difficult to show that every
element of $\scr S$ can be factored into primitive 
elements. This is a
consequence of an old theorem of Bochner [B], though 
Selberg [S] and more
recently Conrey and Ghosh [CG] seem to have found it 
independently.

Selberg conjectures:
\dfn{Conjecture A} For all $F\in\scr S$, there exists a 
positive integer
$n_F$ such that
$$\sum_{p\leq x}\frac{|a_p(F)|^2}{p}=n_F\log\log x+O(1).$$
In Proposition 2.5, we shall describe $n_F$ more 
explicitly.\enddfn
\dfn{Conjecture B} (i) For any primitive function $F$, 
$n_F=1$ so that
$$\sum_{p\leq x}\frac{|a_p(F)|^2}{p}=\log\log x+O(1).$$

(ii) For two distinct primitive functions $F$ and $F'$,
$$\sum_{p\leq x}\frac{a_p(F)\overline{a_p(F')}}{p}=O(1).$$

Thus, in some sense, the primitive functions form an 
orthonormal system.
\enddfn

In his paper [S], Selberg investigates the consequences of 
his conjectures
to the value distribution of $\log F(\sigma+it)$ for 
$\sigma=1/2$ or $\sigma$
very near to $1/2$. Selberg also conjectures the analogue 
of the Riemann
hypothesis for the functions $F\in\scr S$.

It is not difficult to see that Conjecture B implies 
Conjecture A. By
Proposition 2.4, Conjecture B also implies that the 
factorization into
primitives in $\scr S$ is unique. It seems central, 
therefore, to classify
the primitive functions.

To this end, it is natural to define the {\it dimension} 
of $F$ as
$$\dim F=2\alpha_F$$
where
$$\alpha_F=\sum^d_{i=1}\alpha_i.$$
By Proposition 2.2, this concept is well defined. Selberg 
conjectures that the
dimension of $F$ is always a nonnegative integer. This 
question was
previously raised by Vign\'eras [V].

Bochner's work can be used to classify primitive functions 
of dimension
one in the case $\alpha_1=1/2$. They are (after a suitable 
translation) the
classical zeta function of Riemann and the classical 
$L$-functions of
Dirichlet. If the $\alpha_i$ are rational numbers, then 
this is also a
complete list of (primitive) functions of dimension one 
(see [M2]).

The main goal of this paper is to show that Conjecture B 
implies Artin's
conjecture concerning the holomorphy of nonabelian 
$L$-series attached
to irreducible Galois representations. More precisely, let 
$k$ be an
algebraic number field and $K/k$ a finite Galois extension 
with group
$G$. Let $\rho$ be an irreducible representation on the 
$n$-dimensional
complex vector space $V$. For each prime ideal $\germ p$ 
of $k$, let
$V_{\germ p}$
be the subspace of $V$ fixed by the inertia group 
$I_{\germ p}$ of
$\germ p$. Set
$$L_{\germ p}(s,\rho)=\roman{det}(1-\rho(\sigma_{\germ p})
\bold N\germ p^{-s}|V_{\germ p})^{-1}$$
where $\sigma_{\germ p}$ is the Frobenius automorphism of 
the prime ideal
$\germ p$ of $k$
and $\bold N$ is the absolute norm from $k$ to $Q$. Define
$$L(s,\rho;K/k)=\prod_{\germ p}L_{\germ p}(s,\rho).$$
(Sometimes we write $L(s,\rho)$ if the field extension is 
clear. Since
$L(s,\rho)$ depends only on the character $\chi$ of 
$\rho$, we will also
sometimes write $L(s,\chi)$ or $L(s,\chi,K/k)$ for 
$L(s,\rho,K/k)$.)
Clearly, the Artin $L$-function $L(s,\rho)$ is a product 
of $L$-functions
attached to irreducible constituents of $\rho$. We have:
\dfn{Artin's conjecture} If $\rho$ is irreducible 
$\not=1$, then $L(s,\rho,
K/k)$ extends to an entire function of $s$.\enddfn

This embodies one of the central problems of number 
theory. Langlands [L1] has
enunciated an even stronger conjecture: if $n=\dim\rho$, 
then there is a
cuspidal automorphic representation $\pi$ of $GL_n(A_k)$, 
where
$A_k$ is the adele ring of $k$, so that
$$L(s,\rho,K/k)=L(s,\pi).$$
This is called the {\it Langlands reciprocity conjecture}. 
In case $\rho$
is one-dimensional, the conjecture reduces to the 
celebrated Artin
reciprocity law.

Using Artin's reciprocity law, we will show that 
Conjecture B implies even the
Langlands reciprocity conjecture if $K/Q$ is solvable. The 
constraint of
solvability arises from the work of Arthur and Clozel 
[AC], who showed that
the maps of base change and automorphic inductions for 
automorphic
representations exist if $K/k$ is cyclic of prime degree. 
The maps are,
however, conjectured to exist in general. If such is the 
case, then the
Langlands reciprocity conjecture is true in its full 
generality (see [M1]).
In fact, it will be apparent from Theorem 4.3 and Remark 
4.4 that if the
Dedekind zeta function of $K$ is the $L$-function of an 
automorphic
representation over $Q$, then Conjecture B implies the 
Langlands reciprocity
conjecture.

The concept of a primitive function is an important one. 
The Selberg
conjectures imply that $L$-functions attached to 
irreducible cuspidal
automorphic representations over $Q$ are primitive. It 
would be 
extremely interesting to determine whether there are 
elements of $\scr S$
which do not arise from automorphic representations.
\ah 2. Preliminaries\endah

We record in this section the results of Bochner [B], 
Selberg [S], and
Conrey-Ghosh [CG].
\thm\nofrills{Proposition 2.1\ {\rm(Bochner [B])}}. If 
$F\in\scr S$ and
$\alpha_F>0$, then $\alpha_F\geq 1/2$.\ethm
\rem{Remark} In their paper, Conrey and Ghosh [CG] give a 
simple proof
of this and also treat the case $\alpha_F=0$. They prove 
that the constraint
$b_n=O(n^\theta)$ for some $\theta<1/2$ implies there is 
no element in $\scr S$
with $\alpha_F=0$ except the constant function 1.\endrem
\thm\nofrills{Proposition 2.2\ {\rm(Selberg [S])}}. Let 
$N_F(T)$ be the number
of zeros $\rho=\beta+i\gamma$ of $F(s)$ satisfying 
$0<\gamma\leq T$. Then
$$N_F(T)=\frac{\alpha_F}{\pi}T(\log T+c)+S_F(T)+O(1),$$
where $c$ is a constant and $S_F(T)=O(\log T)$.\ethm

If $F=F_1F_2$, then clearly $N_F(T)=N_{F_1}(T)+N_{F_2}(T)$ 
so that
$\alpha_F=\alpha
_{F_1}+\alpha_{F_2}$. Thus, if $F$ is such that 
$\alpha_F<1$, then $F$ is
necessarily primitive. The following is now immediate.
\thm\nofrills{Proposition 2.3\ {\rm[CG]}}. Every $F\in\scr 
S$
has a factorization into primitive functions.\ethm
\demo{Proof} If $F$ is not primitive, then $F=F_1F_2$, 
and, by the above,
$\alpha_F=\alpha_{F_1}+\alpha_{F_2}$. By Proposition 2.1, 
each of
$\alpha_{F_1}$, $\alpha_{F_2}$ is strictly less than 
$\alpha_F$. Continuing
this process, we find that the process terminates because 
Proposition 2.1
implies the number of factors is $\leq 2\alpha_F$. This 
completes the proof.

We see immediately that the Riemann zeta function and the 
classical
Dirichlet functions $L(s,\chi)$ with $\chi$ a primitive 
character are
primitive in the sense of Selberg. Indeed, the 
$\Gamma$-factor appearing
in the functional equation is $\Gamma(s/2)$ or $\Gamma((s+
1)/2)$, and the
result is now clear from Proposition 2.1.

Conjecture B forces the factorization in Proposition 2.3 
to be unique.
Indeed, suppose that $F$ had two factorizations into 
primitive functions:
$$F=F_1\cdots F_r=G_1\cdots G_t$$
where $F_1,\dotsc,F_r$, $G_1,\dotsc,G_t$ are primitive 
functions. Without
loss of generality, we may suppose that no $G_i$ is an 
$F_1$. But then,
$$a_p(F_1)+\cdots+a_p(F_r)=a_p(G_1)+\cdots+a_p(G_t)$$
so that
$$\sum_{p\leq x}\frac{\overline{a_p(F_1)}(a_p(F_1)+\cdots+
a_p(F_r))}
{p}=\sum_{p\leq x}\frac{\overline{a_p(F_1)}(a_p(G_1)+
\cdots+a_p(G_t))}
{p}.$$
As $x\rightarrow\infty$, the left-hand side tends to 
infinity, whereas the
right-hand side is bounded, since no $G_i$ is an $F_1$. 
This contradiction
proves:\enddemo
\thm\nofrills{Proposition 2.4\ {\rm[CG]}}. Conjecture $B$
implies that every $F\in\scr S$ has a unique factorization 
into
primitive functions.\ethm

In the next proposition, we describe $n_F$.
\thm{Proposition 2.5} $(\roman a)$ If $F\in\scr S$ and 
$F=F^{e_1}
_1\cdots F_r^{e_r}$ is a factorization into primitive 
functions, then
Conjecture $B$ implies $n_F=e^2_1+\cdots+e^2_r$.

$(\roman b)$ Conjecture $B$ implies that $F$ is primitive 
if and only if
$n_F=1$.\ethm
\demo{Proof} We have
$$a_p(F)=\sum^r_{i=1}e_ia_p(F_i),$$
so computing the asymptotic behavior
of $\sum_{p\leq x}|a_p(F)|^2/p$ and using Conjecture B 
yield the result.
\enddemo
\ah 3. Artin's conjecture\endah

We begin by showing that Selberg's conjectures imply the 
holomorphy
of nonabelian $L$-functions. In the notation of \S1, let 
$\chi$ be the
character of the representation $\rho$. We will write 
$L(s,\chi,K/k)$
for $L(s,\rho,K/k)$.
\thm{Theorem 3.1} Conjecture $B$ implies Artin's 
conjecture.\ethm
\demo{Proof} We adhere to the notation introduced in \S1. 
Let $\widetilde{K}$
be the normal closure of $K$ over $Q$. Then, 
$\widetilde{K}/k$ is Galois,
as well as $\widetilde{K}/Q$, and $\chi$ can be thought of 
as a character
$\widetilde{\chi}$ of $\roman{Gal}(\widetilde{K}/k)$. By 
the property
of Artin $L$-functions (see [A]),
$$L(s,\widetilde{\chi},\widetilde{K}/k)=L(s,\chi,K/k).$$
Moreover, if $\roman{Ind}\,\widetilde{\chi}$ denotes the 
induction of
$\widetilde{\chi}$ from $\roman{Gal}(\widetilde{K}/k)$ to 
$\roman{Gal}
(\widetilde{K}/Q)$, then
$$L(s,\widetilde{\chi},\widetilde{K}/k)=L(s,\roman{Ind}\,%
\widetilde{\chi},
\widetilde{K}/Q),$$
by the invariance of Artin $L$-functions under induction. 
Hence, we can write
$$L(s,\chi,K/k)=\prod_\varphi L(s,\varphi,\widetilde{K}/Q)
^{m(\varphi)}$$
where the product is over irreducible characters $\varphi$ 
of $\roman{Gal}
(\widetilde{K}/Q)$ and $m(\varphi)$ are positive integers. 
To prove
Artin's conjecture, it suffices to show that 
$L(s,\varphi,\widetilde{K}/Q)$
is entire for each irreducible character $\varphi$ of 
$\roman{Gal}(\widetilde{K}/
Q)$. By Brauer's induction theorem and the Artin 
reciprocity law, we can
write
$$L(s,\varphi,\widetilde{K}/Q)=\frac{L(s,\chi_1)}{L(s,%
\chi_2)}$$
where $\chi_1$ and $\chi_2$ are characters of 
$\roman{Gal}(\widetilde{K}/Q)$
and $L(s,\chi_1),L(s,\chi_2)$ are entire functions, being 
products
of Hecke $L$-functions. Thus, they belong to $\scr S$ and, 
hence, by
Proposition 2.4 have a unique factorization into primitive 
functions.
We can therefore write
$$L(s,\varphi)=\prod^m_{i=1}F_i(s)^{e_i},\qquad e_i\in Z.$$
By comparing the $p\roman{th}$ Dirichlet coefficient of 
both sides, we get
$$\varphi(p)=\sum^m_{i=1}e_ia_p(F_i),$$
from which we obtain
$$\sum_{p\leq x}\frac{|\varphi(p)|^2}{p}=\sum_{p\leq 
x}\frac{1}{p}
\left|\sum^m_{i=1}e_ia_p(F_i)\right|^2.$$
Conjecture B gives the asymptotic behaviour of the 
right-hand side:
$$\sum_{p\leq 
x}\frac{|\varphi(p)|^2}{p}=\left(\sum^m_{i=1}e^2_i\right)
\log\log x+O(1).$$
Decompose the sum on the left-hand side according to the 
conjugacy class
$C$ of $\roman{Gal}(\widetilde{K}/Q)$ to which the 
Frobenius automorphism
$\sigma_p$ belongs:
$$\sum_{p\leq x}\frac{|\varphi(p)|^2}{p}=\sum_C
|\varphi(g_C)|^2\sum\Sb p\leq x\\
\sigma_p\in C\endSb\frac{1}{p},$$
where $g_C$ is any element of $C$. By the Chebotarev 
density theorem
$$\sum\Sb p\leq x\\
\sigma_p\in C\endSb\frac{1}{p}=\frac{|C|}{|G|}\log\log x+
O(1).$$
Hence,
$$\sum_{p\leq 
x}\frac{|\varphi(p)|^2}{p}=\sum_C\frac{|C|}{|G|}
|\varphi(g_C)|^2\log\log x+O(1).$$
But $\varphi$ is irreducible, so
$$\sum_C\frac{|C|}{|G|}|\varphi(g_C)|^2=(\varphi,%
\varphi)=1.$$
Therefore, the left-hand side is
$$\log\log x+O(1)$$
as $x\rightarrow\infty$. We deduce that 
$\sum^m_{i=1}e^2_i=1$,
from which follows $m=1$ and $e_1=\pm 1$. Thus, 
$L(s,\varphi)=F(s)$ or
$1/F(s)$, where $F(s)$ is primitive and analytic 
everywhere except possibly
at $s=1$. However, $L(s,\varphi)$ has trivial zeros, and 
so the latter possibility
cannot arise. We conclude that $L(s,\varphi)=F(s)$ is 
primitive and
entire.\enddemo
\thm{Corollary 3.2} Let $K/Q$ be Galois, and let $\chi$ be 
an irreducible
character of $\roman{Gal}(K/Q)$. Conjecture B implies that
$L(s,\chi)$ is primitive.\ethm
\demo{Proof} This is evident from the last line in the 
proof of the previous
lemma. Or we can derive it as follows. By the previous 
theorem, $F=L(s,\chi)
\in\scr S$; and by the Chebotarev density theorem, 
$n_F=1$. The result
now follows from Proposition 2.5(b).

Of course, Dedekind's conjecture that the zeta function of 
a number field is
always divisible by the Riemann zeta function follows from 
Artin's
conjecture. However, it is rather interesting to note that 
the unique
factorization conjecture is sufficient to deduce this. 
Indeed, if $K$ is
a number field, $\widetilde{K}$ its Galois closure, and 
$\zeta_K(s)$ is the
Dedekind zeta function of $K$, then 
$\zeta_{\widetilde{K}}(s)/\zeta_K(s)=
F(s)$ is entire by the Aramata--Brauer theorem. By the 
same theorem,
$\zeta_{\widetilde{K}}(s)/\zeta(s)=G(s)$ is also entire. 
Since $\zeta(s)$
is primitive, it appears as a primitive factor in 
$\zeta_{\widetilde{K}}
(s)=\zeta(s)G(s)$. Since 
$\zeta_{\widetilde{K}}(s)=\zeta_K(s)F(s)$ and $F$
is entire, $\zeta(s)$ must appear in the unique 
factorization of $\zeta_K(s)$.
This is Dedekind's conjecture.\enddemo
\ah 4. Langlands reciprocity conjecture\endah

We begin by giving a very brief description of 
$L$-functions attached to
automorphic representations of $GL_n$. Such a description 
is bound to be
incomplete, so we refer the reader to [Co] for details. 
Let $n\geq 1$,
$G=GL_n$, and $A$ be any commutative ring with identity. 
We denote by
$G(A)$ the group of $n\times n$ matrices over $A$ whose 
determinant is a
unit in $A$. Denote by $M_r$ a generic $r\times r$ matrix 
and by $I_r$
the $r\times r$ identity matrix. The standard parabolic 
subgroups of
$GL_n$ are in one-to-one correspondence with the 
partitions of $n=n_1
+\cdots+n_r$. The standard parabolic subgroup 
corresponding to a partition
$n=n_1+\cdots+n_r$ consists of matrices of the form
$$\pmatrix M_{n_1}&*&*\\
&\ddots&*\\
&&M_{n_r}\endpmatrix,$$
and any parabolic subgroup is a $GL_n$ conjugate of a 
standard parabolic
subgroup. Any parabolic subgroup $P$ has a decomposition 
(called the
Levi decomposition) of the form $P=MN$ where $N$ is the 
unipotent radical
of $P$. $M$ is called the Levi component of $P$. In the 
case of the standard
parabolic, $M$ and $N$ can be described as consisting of 
matrices of the form
$$\pmatrix M_{n_1}&&\\
&\ddots&\\
&&M_{n_r}\endpmatrix\qquad\pmatrix I_{n_1}&*&*\\
&\ddots&*\\
&&I_{n_r}\endpmatrix$$
respectively. We will write $N_P$ to denote the unipotent 
radical in the
Levi decomposition of a parabolic subgroup $P$.

If $k$ is an algebraic number field, then the adele ring 
of $k$, denoted
$A_k$, is a commutative ring with identity. It is defined 
in the following
way. For each place $v$, let $k_v$ be the completion of 
$k$ at $v$. As a set
$A_k$ consists of all infinite tuples $(x_v)$ where $v$ 
ranges over all
places of $k$ with $x_v$ in $k_v$ and $x_v$ in the ring of 
$v$-adic integers,
$O_v$ for all but finitely many places. One defines 
addition and multiplication
componentwise which makes it into a commutative ring. We 
impose the adelic
topology by declaring for each finite set $S$ of places 
containing the
archimedean places,
$$\prod_{v\in S}k_v\times\prod_{v\not\in S}O_v,$$
with the product topology, as a basic neighborhood of the 
identity. This
makes $A_k$ into a locally compact topological ring. One 
can think of $k$ as
embedded in the ring $A_k$ via the map
$$x\mapsto(x,x,\dotsc).$$
The adelic topology on $G(A_k)$ is similarly defined. For 
each finite set
$S$ as above, we declare
$$\prod_{v\in S}G(k_v)\times\prod_{v\not\in S}G(O_v),$$
with the product topology, to be a basic neighborhood of 
the identity. The
adelic topology on $G(A_k)$ makes it into a locally 
compact group in which
$G(k)$, embedded diagonally, is a discrete subgroup of 
$G(A_k)$. The coset
space, $G(k)\backslash G(A_k)$, with quotient topology 
does not have finite
volume with respect to any $G(A_k)$ invariant measure. To 
rectify this, define
$$Z=\left\{\pmatrix
z&&\\
&\ddots&\\
&&z\endpmatrix\colon\ z\in A^\times_k
\right\}.$$
Then, the quotient $ZG(k)\backslash G(A_k)$ has finite 
volume with respect to
any $G(A_k)$ invariant measure.

A character of $k^\times\backslash GL_1(A_k)$ is called a 
Grossencharacter.
Now fix a Grossencharacter $\omega$ of $k$. We can 
consider the Hilbert
space
$$L^2(G(k)\backslash G(A_k),\omega)$$
of measurable functions $\varphi$ on $G(k)\backslash 
G(A_k)$ satisfying
\roster
\item"(i)" $\varphi(zg)=\omega(z)\varphi(g)$, $z\in Z$, 
$g\in G(k)\backslash G(A_k)$;
\item"(ii)" $\int_{ZG(k)\backslash 
G(A_k)}|\varphi(g)|^2\,dg<\infty$.\newline
The subspace of cusp forms $L^2_0(G(k)\backslash 
G(A_k),\omega)$ is defined
by the extra condition
\item"(iii)" for all parabolic subgroups $P$ of $G(A_k)$,
$$\int_{N_P(k)\backslash 
N_P(A_k)}\varphi(ng)\,dn=0$$\endroster
for every $g\in G(A_k)$.

Let $R$ be the right regular representation of $G(A_k)$ on 
$L^2(G(k)\backslash
G(A_k),\omega)$. Thus,
$$(R(g)\varphi)(x)=\varphi(xg)$$
for $\varphi\in L^2(G(k)\backslash G(A_k,\omega)$ and 
$x,g\in G(A_k)$. This is a
unitary representation of $G(A_k)$. An automorphic 
representation is a
subquotient of the right regular representation of 
$G(A_k)$ on $L^2(G(k)
\backslash G(A_k),\omega)$. A cuspidal automorphic 
representation is a
subrepresentation of the right regular representation of 
$G(A_k)$ on
$L^2_0(G(k)\backslash G(A_k),\omega)$.

A representation of $G(A_k)$ is called admissible if its 
restriction to the
maximal compact subgroup
$$K=\prod_{v\,\roman{complex}}U(n,C)\times\prod_{v\,%
\roman{real}}O(n,R)
\times\prod_{v\,\roman{finite}}GL_n(O_v)$$
contains each irreducible representation of $K$ with only 
finite
multiplicity.

To understand the structure of these representations, we 
need the notion
of a restricted tensor product of representations 
originally introduced
in [JL]. Let $\{W_v|v\in V\}$ be a family of vector spaces 
indexed by the
set $V$. Let $V_0$ be a finite subset of $V$. For each 
$v\in V\backslash V_0$,
let $x_v$ be a nonzero vector in $W_v$. For each finite 
subset $S$ of $V$
containing $V_0$, let $W_s=\bigotimes_{v\in S}W_v$. If 
$S\subset S'$,
let $f_{S,S'}\colon\ W_S\rightarrow W_{S'}$ be defined by
$$\bigotimes_{v\in S}w_v\mapsto\bigotimes_{v\in 
S}w_v\bigotimes_{v\in
S'\backslash S}x_v.$$
Then, the restricted tensor product of the $W_v$ with 
respect to the $x_v$,
$W=\bigotimes_{x_v}\!W_v$, is defined to be the inductive 
limit $W=\ind\lim
_SW_S$. It is known [F, p. 181] that any irreducible, 
admissible representation
of $GL_n(A_k)$ can be written as a restricted tensor 
product $\bigotimes
_v\pi_v$ where $\pi_v$ is an irreducible representation of 
$GL_n(k_v)$.
Moreover, the factors $\pi_v$ are unique up to equivalence.

For each finite $v$, one can construct certain 
representations of
$GL_n(k_v)$ in the following way. One has the Borel subgroup
$$B(k_v)=\left\{b=\pmatrix b_1&*&*\\
&\ddots&*\\
&&b_n\endpmatrix\right\}\subseteq GL_n(k_v)$$
if $GL_n(k_v)$. For any $n$-tuple, $z=(z_1,\dotsc,z_n)\in 
C^n$, define
$\chi_z$ on $B(k_v)$ by
$$\chi_z(b)=|b_1|^{z_1}_v\cdots|b_n|^{z_n}_v,$$
where $|\boldcdot|_v$ denotes the $v$-adic metric. This 
gives a
quasi-character on $B(k_v)$. Let $\widetilde{\pi}_{v,z}$ 
be the
representation of $GL_n(k_v)$ obtained by inducing 
$\chi_z$ from
$B(k_v)$ to $GL_n(k_v)$. Recalling the definition of an 
induced representation,
we see that this means the following: it is the space of 
locally constant
functions $\varphi$ on $GL_n(k_v)$ such that for $b\in 
B(k_v)$,
$x\in GL_n(k_v)$,
$$\varphi(bx)=\chi_z(b)\left(\prod^n_{i=1}|b_i|_v^{(n+
1)/2-i}\right)\varphi(x),$$
and $\widetilde{\pi}_{v,z}$ acts on $\varphi$ via
$$(\widetilde{\pi}_{v,z}(y)\varphi)(x)=\varphi(xy).$$
We shall assume that
$$\roman{Re}(z_1)\geq\roman{Re}(z_2)\geq\cdots\geq%
\roman{Re}(z_n).$$
Then, a special case of the Langlands classification [BW, 
Chapter 11, \S2]
shows that $\widetilde{\pi}_{v,z}$ has a unique 
irreducible quotient
$\pi_{v,z}$. The representations obtained in this way are 
called the unramified
principal series. They are precisely the representations 
of $GL_n(k_v)$ whose
restriction to $GL_n(O_v)$ contains the trivial 
representation. If $\pi_v$
is any representation of $GL_n(k_v)$ equivalent to 
$\pi_{v,z}$, we can
define
$$A_v=\roman{diag}(Nv^{-z_1},\dotsc,Nv^{-z_n})\in GL_n(C)$$
where $Nv$ denotes the norm of $v$. One can show that 
$A_v$ depends only
on the equivalence class of $\pi_v$.

Now let $\pi=\bigotimes_v\pi_v$ be an irreducible, 
admissible, automorphic
representation of $GL_n(A_k)$. Since $\pi$ is admissible, 
there is a finite
set $S$ of places (containing the infinite places) such 
that $\pi_v$ belongs
to the unramified principal series for $v\not\in S$. For 
such $v$, define
$$L_v(s,\pi)=\det(1-A_vNv^{-s})^{-1}.$$
Set
$$L_S(s,\pi)=\prod_{v\not\in S}L_v(s,\pi).$$
Using the theory of Eisenstein series, Langlands [L3] 
established the
meromorphic continuation of $L_S(s,\pi)$. It is possible 
to define
$L_v(s,\pi_v)$ for $v\in S$ such that the complete 
$L$-function
$$L(s,\pi)=\prod_vL_v(s,\pi_v)$$
has a meromorphic continuation and functional equation. 
Moreover, if $\pi$
is cuspidal, then $L(s,\pi)$ extends to an entire function 
unless $n=1$
and $\pi$ is of the form $|\boldcdot|^t$ for some $t\in 
C$. This was proved
by Godement and Jacquet [GJ]. The {\it Ramanujan 
conjecture} asserts that
the eigenvalues of $A_v$ are of absolute value 1 for 
cuspidal automorphic
representations. This is still one of the main open 
problems in the theory
of automorphic representations. If the Ramanujan 
conjecture is true, then these
$L$-functions belong to the Selberg class.

For $n=1$ this is the classical theory of Hecke's 
$L$-series attached to
Grossencharacters, and Tate's thesis offers the important 
perspective of an
adelic reformulation of Hecke's work. For $n=2$ and $K=Q$, 
the analytic
continuation and functional equation of the corresponding 
$L$-function
are due to Hecke and Maass. The generalization to global 
fields and the
adelic formulation is the work of Jacquet and Langlands 
[JL]. The general
case for $GL_n$ is due to Godement and Jacquet [GJ]. For 
$n=3$, a simple
proof can be found in the work of Hoffstein and Murty 
[HM]. A gentle
introduction to the Langlands program can be found in [M1].

In this adelic framework, Langlands conjectures that every 
Artin
$L$-series $L(s,\rho)$ is an $L(s,\pi)$ for some 
automorphic representation
$\pi$ on $GL_n(A_k)$ where $n=\dim\rho$. Thus, to each 
$\rho$ there should
be an automorphic representation $\pi(\rho)$. In case 
$\deg\rho=1$, this is
Artin's reciprocity law, since Hecke's Grossencharacters 
are automorphic
representations of $GL_1(A_k)$.

Let $K/k$ be a Galois extension, and let 
$G=\roman{Gal}(K/k)$. Let $H$
be a subgroup of $G$ and $K^H$ its fixed field. If $\psi$ 
is a representation
of $H$, then $L(s,\psi,K/K^H)$ is the Artin $L$-series 
belonging to the
extension $K/K^H$. We have already noted that Artin 
$L$-series are invariant
under induction:
$$L(s,\roman{Ind}^G_H\psi,K/k)=L(s,\psi,K/K^H).$$
If the Langlands reciprocity conjecture is true, this last 
property should
imply a corresponding property for $L$-series attached to 
automorphic
representations. This leads to the concept of the base 
change $L$-function
which we now define.

To this end, let us first describe how the map 
$\rho\mapsto\pi(\rho)$ behaves
under restriction to a subgroup. If we denote by $\scr 
A(GL_n(A_k))$ the
space of automorphic representations of $GL_n(A_k)$, then 
by the
reciprocity conjecture, we should have
$$\rho|_H\mapsto\pi(\rho|_H)\in\scr A(GL_n(A_{K^H})).$$
What is $L(s,\rho|_H,K/K^H)$? Since
$$\roman{Ind}^G_H(\rho|_H\otimes\psi)=\rho\otimes%
\roman{Ind}^G_H\psi,$$
we find
$$L(s,\roman{Ind}^G_H(\rho|_H\otimes\psi),K/k)=L(s,\rho%
\otimes
\roman{Ind}^G_H\psi,K/k).$$
But the former $L$-function is 
$L(s,\rho|_H\otimes\psi,K/K^H)$ by the
invariance of Artin $L$-series. Therefore,
$$L(s,\rho|_H,K/K^H)=L(s,\rho\otimes\roman{Ind}^G_H1,K/k).$$
But $\roman{Ind}^G_H1=\roman{reg}_{G/H}$ is the 
permutation representation
on the cosets of $H$. This suggests that we make the 
following definition.
Let $\pi\in\scr A(GL_n(A_k))$. Define, for unramified $v$,
$$L_v(s,B(\pi))=\det(1-A_v\otimes\deg_{G/H}(%
\sigma_v)Nv^{-s})^{-1},$$
where $\sigma_v$ is the Artin symbol of $v$. For ramified 
$v$, see
[L2] or [AC] for the definition of $L_v(s,B(\pi))$. Set
$$L(s,B(\pi))=\prod_vL_v(s,B(\pi)).$$
{\it B}$(\pi)$ should correspond to an element of $\scr 
A(GL_n(A_M))$
where $M=K^H$. The problem of base change is to determine 
when this map
exists.

For $n=2$ and $M/k$ cyclic, this is a theorem of Langlands 
[L2]. For
arbitrary $n$, $M/k$ cyclic, it is the main theorem of 
Arthur and Clozel
[AC].

The reciprocity conjecture forces another conjectural 
property for the
automorphic $L$-functions. Suppose that $\psi$ is a 
representation of
$H$. Corresponding to $\psi$ there should be a $\pi\in\scr 
A(GL_n(A_M))$
where $n=\roman{dim}\,\psi$. But the invariance of Artin 
$L$-series under
induction implies that there should be an $I(\pi)\in\scr 
A(GL_{nr}(A_k))$
$(r=[G:H])$ so that
$$L(s,I(\pi))=L(s,\roman{Ind}^G_H\psi,K/k).$$
This map $\pi\mapsto I(\pi)$, called the automorphic 
induction map,
is conjectured to exist. Again, Arthur and Clozel [AC] 
showed this exists
when $M/k$ is cyclic and arbitrary $n$. Thus, if $M/k$ is 
contained
in a solvable extension of $k$, the base change maps and 
automorphic
induction maps exist. We can summarize this discussion in 
the following.
\thm\nofrills{Lemma 4.1 {\rm[AC]}}. Let $E/F$ be a cyclic 
extension
of prime degree $p$ of algebraic number fields, and let 
$\pi$ be an
automorphic representation of $GL_n(A_F)$. Then there is 
the base change
lift $B(\pi)$ of $\pi$ such that $B(\pi)$ is an 
automorphic representation
of $GL_n(A_E)$. If $\sigma$ is an automorphic 
representation of $GL_n(A_E)$,
then the automorphic induction $I(\sigma)$ exists.\ethm
We now have the requisite background to prove:
\thm{Lemma 4.2} Assume Conjecture $B$. If $\pi$ is an 
irreducible cuspidal
automorphic representation of $GL_n(A_Q)$ that satisfies 
the Ramanujan
conjecture, then $L(s,\pi)$ is primitive.\ethm
\demo{Proof} Let $\pi_1$ and $\pi_2$ be two irreducible 
cuspidal automorphic
representations of $GL_n(A_Q)$. We recall the work of 
Jacquet and Shalika
[JS1, JS2] and Shahidi [Sh]. Let $S$ be a finite set of 
primes such that
$\pi_1,\pi_2$ are unramified outside $S$. We form the 
$L$-function
$$L_S(s,\pi_1\otimes\pi_2)=\prod_{v\not\in 
S}\det(1-A_{v,1}\otimes A
_{v,2}Nv^{-s})^{-1},$$
where $A_{v,1}$ and $A_{v,2}$ are the diagonal matrices 
associated to
$\pi_1$ and $\pi_2$ respectively, as indicated above. 
Then, we will use
the following properties, derived in [JS1, JS2]:

(a) The Euler product $L_S$ is absolutely convergent for 
$\roman{Re}\,s>1$.

(b) Let $X$ be the set of $s$ on the line 
$\roman{Re}\,s=1$ such that
$\pi_1\otimes|\boldcdot|^{s-1}$ is equivalent to 
$\widetilde{\pi}_2$,
the contragredient of $\pi_2$. (Thus, $X$ contains at most 
one point.)
Then the function $L_S$ extends continuously to the line 
$\roman{Re}\,s=1$
with $X$ removed. Moreover, it does not vanish there.

(c) If $s_0\in X$, the limit
$$\lim\Sb s\rightarrow s_0\\
\roman{Re}\,s\geq 1\endSb(s-s_0)L_S(s,\pi_1\otimes\pi_2)$$
exists and is finite and nonzero.

These facts are proved in [JS1, Theorem 5.3, JS2, 
Proposition 3.6], and the
nonvanishing results are proved by Shahidi [Sh]. When 
$\pi_1$ and
$\widetilde{\pi}_2$ are distinct, this means that 
$L_S(s,\pi_1\otimes
\pi_2)$ is analytic and nonvanishing for 
$\roman{Re}\,s=1$, whereas
$L_S(s,\pi_1\otimes\widetilde{\pi}_1)$ has a simple pole 
at $s=1$ and is
nonvanishing for $\roman{Re}\,s=1$. In analytic terms, 
this implies
$$\sum_{p\leq x}\frac{|a_p(\pi_i)|^2}{p}=\log\log x+O(1)$$
and
$$\sum_{p\leq 
x}\frac{a_p(\pi_1)\overline{a_p(\pi_2)}}{p}=O(1).$$
By Proposition 2.5(b), we deduce the desired result.\enddemo

We are now ready to prove the main theorem of this section.
\thm{Theorem 4.3} Assume Conjecture B. Let $K$ be a Galois 
extension
of $Q$ with solvable group $G$, and let $\chi$ be an 
irreducible character
of $G$ of degree $n$. Then there is an irreducible 
cuspidal automorphic
representation $\pi$ of $GL_n(A_Q)$ such that 
$L(s,\chi)=L(s,\pi)$.\ethm
\demo{Proof} The Dedekind zeta function of $K$ is a Hecke 
$L$-function.
In fact, it is an automorphic $L$-function of $GL_1(A_K)$. 
Since $K$
is solvable over $Q$, there is a chain of cyclic extensions,
$$Q=K_0\subset K_1\cdots\subset K_{m-1}\subset K_m=K$$
so that $[K_i\colon\ K_{i-1}]$ has prime degree. By Lemma 
4.1 and automorphic
induction, we can identify $\zeta_K(s)$ as an automorphic 
$L$-function,
first of $K_{m-1}$ and then successively of $K_i$ until we 
deduce
that it is an automorphic $L$-function of $Q$. Thus, as in 
[AC], we have
the factorization
$$\zeta_K(s)=\prod^r_{i=1}L(s,\pi_i)^{e_i}$$
where $\pi_i$ are irreducible cuspidal automorphic 
representations of
appropriate degree over $Q$ and $e_i$ are positive 
integers. Since the
eigenvalues of the $p$-Euler factors of the left-hand side 
are of
absolute value 1, so are the eigenvalues of the $p$-Euler
factors on the right-hand side. Thus, each $L(s,\pi_i)$ 
satisfies the
Ramanujan conjecture. By Lemma 4.2, each factor is 
primitive. Thus, the
above is the unique factorization of $\zeta_K(s)$ into 
primitive factors.
On the other hand, we have the classical factorization
$$\zeta_K(s)=\prod_\chi L(s,\chi)^{\chi(1)},$$
where the product is over the irreducible characters of 
$G$. By Theorem 3.1
and Corollary 3.2, each $L(s,\chi)$ is in $\scr S$ and, 
moreover, primitive.
Thus, each $L(s,\chi)$ must be an $L(s,\pi)$ by uniqueness 
of factorization.
This is the reciprocity law.\enddemo
\rem{Remark {\rm4.4}} If $K$ were an arbitrary Galois 
extension of $Q$ and we knew
that $\zeta_K(s)$ was automorphic, the reciprocity 
conjecture would follow
in the same way on the assumption of Conjecture B.\endrem
\ah 5. Concluding remarks\endah

The Selberg conjectures refer to the analytic behavior of 
Dirichlet series
at the edge of the critical strip. There are other 
conjectures relating special
values of Dirichlet series inside the critical strip, 
namely, the Deligne
conjectures and the Birch--Swinnerton-Dyer conjectures to 
cite specific
instances. A consequence of Conjecture B is that if $F$ is 
any primitive
function which is not the Riemann zeta function, then
$$\sum_{p\leq x}\frac{a_p(F)}{p^{1+it}}=O(1).$$
In particular, no primitive function should vanish on 
$\sigma=1$.

Many of our interesting consequences, notably the Artin 
conjecture and the
Langlands reciprocity conjecture, utilized the unique 
factorization conjecture.
Perhaps this can be attacked by other means. Indeed, given 
$r$ distinct
primitive functions $F_1,\dotsc,F_r$, one would expect the 
existence of
complex numbers $s_1,\dotsc,s_r$ such that $F_i(s_j)=0$ if 
and only if
$i=j$. If this were the case, then clearly, the unique 
factorization
conjecture would be true.

The classification of primitive functions is a fundamental 
problem. From the
work of Bochner and Vign\'eras, it follows that if $F$ has 
dimension 1
and all the $\alpha_i$ are rational numbers, then $d=1$ 
and $\alpha_1=1/2$.
It then follows, essentially from the same works, that $F$ 
must either
be the Riemann zeta function or a purely imaginary 
translate of a classical
Dirichlet $L$-function attached to a nontrivial primitive 
character. It is
shown in [M2] that if $\pi$ is an irreducible cuspidal 
automorphic
representation of $GL_2(A_Q)$, then $L(s,\pi)$ is 
primitive if the
Ramanujan conjecture is true. In particular, the 
$L$-function attached
to a normalized holomorphic cuspidal Hecke eigenform is a 
primitive
function which is in Selberg's class (by Deligne's theorem).
\ah Acknowledgment\endah

I would like to thank Professors Barry Mazur, J.-P. Serre, 
and Dipendra
Prasad for their critical comments on an earlier version 
of this paper.

\enddocument